\documentclass[letterpaper]{article} 
\usepackage[draft]{aaai2027} 
\usepackage[hyphens]{url} 
\usepackage{graphicx} 
\usepackage{natbib} 
\usepackage{caption} 
\usepackage{amsmath,amssymb,amsthm,mathtools,bm,pifont}
\usepackage{array}
\usepackage{booktabs}
\usepackage{enumitem}
\usepackage{float}
\usepackage{multirow}
\usepackage{needspace}
\usepackage{xcolor}
\usepackage{cancel}
\usepackage{tabularx}

\providecommand{\fro}{\mathsf{F}}
\providecommand{\ts}{\textstyle}
\providecommand{\A}{\mathbf{A}}
\providecommand{\B}{\mathbf{B}}
\providecommand{\Y}{\mathbf{Y}}
\providecommand{\G}{\mathbf{G}}
\providecommand{\I}{\mathbf{I}}
\providecommand{\M}{\mathbf{M}}
\providecommand{\R}{\mathbf{R}}
\providecommand{\X}{\mathbf{X}}

\providecommand{\hG}{\hat{G}}
\providecommand{\EEE}{\mathbb{E}}
\providecommand{\GGG}{\mathbb{G}}
\providecommand{\MMM}{\mathbb{M}}
 \def\zero{\mathbf{0}}
\providecommand{\varepsilons}{\boldsymbol{\varepsilon}}
\newcommand{\beq}{\begin{eqnarray}}
\newcommand{\eeq}{\end{eqnarray}}
\newcommand{\noi}{\noindent}
\newcommand{\nn}{\nonumber}
\newcommand{\la}{\langle}
\newcommand{\ra}{\rangle}
\def\Rn{\mathbb{R}} 
\def\ts{\textstyle}
\newcommand{\step}[1]{\text{\ding{\numexpr#1+171\relax}}}
\newcommand{\xxxx}[2]{\if\relax\detokenize\expandafter{\romannumeral-`\q#1}\relax#2\else \overset{\step{#1}}{#2}\fi}

\newtheoremstyle{NABLA-plain}
  {\topsep}{\topsep}
  {\itshape}{}
  {\bfseries}{.}
  {0.5em}
  {\thmname{#1}\thmnumber{~#2}\thmnote{~(#3)}}
\theoremstyle{NABLA-plain}
\newtheorem{theorem}{Theorem}
\newtheorem{lemma}{Lemma}
\newtheorem{assumption}{Assumption}
\newcounter{optalg}
\newsavebox{\optalgbox}
\providecommand{\normalem}{}
\newenvironment{keywords}{\par\smallskip\noindent\textbf{Keywords: }\ignorespaces}{\par\smallskip}

\title{NABLA: A Noise-Adaptive Blockwise Left--Right Preconditioning Algorithm for Stochastic Optimization}

\author{
    Yanlin Li\textsuperscript{\rm 1},
    Ganzhao Yuan\textsuperscript{\rm 2\ensuremath{\dagger}}
}
\affiliations{
    \textsuperscript{\rm 1}Southern University of Science and Technology\\
    \textsuperscript{\rm 2}Shenzhen University of Advanced Technology\\
}

\begin{document}
\maketitle

\begin{abstract}

We propose NABLA, a noise-adaptive blockwise left--right preconditioning algorithm for mini-batch stochastic nonconvex optimization. NABLA combines an exponential moving average of stochastic gradients with block-diagonal row and column preconditioners constructed from matrix-valued second-moment statistics. Unlike conventional matrix-preconditioned methods that use fixed decay parameters, NABLA adjusts the second-moment averaging coefficient through a closed-loop feedback rule determined by the cumulative stochastic-gradient energy and its running average. The first-moment coefficient remains fixed, allowing the trajectory-dependent adaptation to act directly on the structured preconditioners. The resulting method is horizon-free, requires no prior knowledge of the noise level, and automatically adjusts the memory of its second-moment estimators to the observed optimization dynamics. Under standard lower boundedness, smoothness, unbiasedness, bounded-variance, and spectral-boundedness assumptions, we establish uniform stability and variation bounds for the time-varying block preconditioners and control their interaction with momentum estimation and parameter updates. These results yield the noise-adaptive stationarity guarantee $\widetilde{\mathcal{O}}(T^{-1/2}+\sigma^{1/2}T^{-1/4})$, which automatically reduces to the nearly deterministic $\widetilde{\mathcal{O}}(T^{-1/2})$ rate when the stochastic-gradient noise vanishes. NABLA therefore provides a principled approach to combining closed-loop trajectory feedback with memory-efficient structured matrix preconditioning.

\end{abstract}

\begin{minipage}{\columnwidth}
\end{minipage}

\section{Introduction}

We consider mini-batch stochastic optimization of the possibly nonconvex matrix problem
\begin{equation}\label{eq:main}
\min_{\X\in\mathbb{R}^{m\times n}}
f(\X):=\mathbb{E}_{\xi\sim\mathcal{D}}[f(\X;\xi)].
\end{equation}
Here $\xi$ is an i.i.d. sample from an unknown distribution $\mathcal D$, and the objective is differentiable and bounded below. Our goal is to produce an iterate with a small expected Frobenius norm of the true gradient.

Adaptive gradient methods have become indispensable in modern stochastic optimization. On the one hand, exponential moving average (EMA) mechanisms, as used in Adam-type methods, stabilize stochastic gradients and improve empirical robustness. On the other hand, Shampoo-type methods exploit matrix-valued second-order statistics to precondition gradients from both the left and the right, thereby capturing richer geometry than coordinate-wise scaling. Despite their practical success, the theoretical understanding of adaptive matrix-preconditioned methods remains less complete in the nonconvex stochastic setting, especially when one asks for a guarantee that simultaneously adapts to noise and recovers the fast deterministic behavior when the noise level is zero.

A key difficulty is that standard adaptive stepsizes are often controlled by historical accumulations. Such accumulations are useful for stability, but they may also decay too aggressively and lead to a suboptimal deterministic rate. This issue becomes particularly important for zero-noise optimality: when $\sigma=0$, the stochastic guarantee should reduce to the nearly deterministic rate rather than remain limited by a stochastic-order term. Therefore, the stepsize mechanism should be both conservative under noise and non-vanishing when the observed gradients become small.

\begin{table*}[!t]
\centering
\caption{Comparison of adaptive optimizers for a matrix parameter
$\mathbf W\in\mathbb R^{m\times n}$ with $m\ge n$.
In NABLA$(k_1,k_2)$, $k_1$ and $k_2$ denote the left and right block sizes;
memory includes the first-moment buffer and all preconditioner statistics.}
\label{tab:structured-optimizers}

\begingroup

\newcommand{\TableScale}{0.92}

\small

\setlength{\tabcolsep}{3.5pt}

\renewcommand{\arraystretch}{1.18}

\newcommand{\FamilyWidth}{1.45cm}
\newcommand{\MethodWidth}{4.45cm}
\newcommand{\AdaptiveWidth}{1.75cm}
\newcommand{\StateWidth}{3.25cm}
\newcommand{\MemoryWidth}{3.50cm}
\newcommand{\RateWidth}{4.15cm}

\resizebox{\TableScale\textwidth}{!}{%

\begin{tabular}{@{}
|>{\centering\arraybackslash}m{\FamilyWidth}
|>{\raggedright\arraybackslash}m{\MethodWidth}
|>{\centering\arraybackslash}m{\AdaptiveWidth}
|>{\raggedright\arraybackslash}m{\StateWidth}
|>{\centering\arraybackslash}m{\MemoryWidth}
|>{\raggedright\arraybackslash}m{\RateWidth}
|@{}}

\hline

\textbf{Family}
&
\textbf{Method and Reference}
&
\shortstack[c]{\textbf{Noise-}\\\textbf{Adaptive?}}
&
\shortstack[l]{\textbf{State /}\\\textbf{Preconditioner}}
&
\textbf{Extra Memory}
&
\shortstack[l]{\textbf{Stochastic Convergence}\\
               \textbf{Rate}}
\\

\hline

\multirow[c]{3}{\FamilyWidth}{\centering Diagonal}
& AdaGrad~\citep{duchi2011adagrad}
& --
& diagonal cumulative statistic
& $\mathcal O(mn)$
& -- (convex regret)
\\ \cline{2-6}

& Adam~\citep{kingma2015adam}
& --
& diagonal EMA
& $\mathcal O(mn)$
& -- (online regret)
\\ \cline{2-6}

& OptEMA-V~\citep{yuan2026optema}
& Yes
& diagonal EMA
& $\mathcal O(mn)$
& \shortstack[l]{
    $\widetilde{\mathcal O}(T^{-1/2} + \sigma^{1/2}T^{-1/4})$
  }
\\ \hline

\multirow[c]{3}{\FamilyWidth}{\centering One-sided}
& Muon~\citep{kim2026convergence}
& --
& right polar/Gram factor
& $\mathcal O(mn+n^2)$
& $\mathcal O(\log(T)/\sqrt T)^{\dagger}$
\\ \cline{2-6}

& One-sided Shampoo~\citep{xie2025structured}
& --
& one Gram factor
& $\mathcal O(mn+n^2)$
& -- (convex analysis)
\\ \cline{2-6}

& OptMuon-A~\citep{yuan2026optmuon}
& Yes
& right polar/Gram factor
& $\mathcal O(mn+n^2)$
& \shortstack[l]{
    $\widetilde{\mathcal O}(T^{-1/2} + \sigma^{1/2}T^{-1/4})$
  }
\\ \hline

\multirow[c]{4}{\FamilyWidth}{\centering Two-sided}
& Shampoo~\citep{gupta2018shampoo}
& --
& full left/right Kronecker
& $\mathcal O(mn+m^2+n^2)$
& -- (convex regret)
\\ \cline{2-6}

& KrADagrad~\citep{mei2023kradagrad}
& --
& full left/right Kronecker
& $\mathcal O(mn+m^2+n^2)$
& -- (convex regret)
\\ \cline{2-6}

& SOAP~\citep{vyas2025soap}
& --
& Shampoo eigenbasis + diagonal EMA
& $\mathcal O(mn+m^2+n^2)$
& -- (empirical study)
\\ \cline{2-6}

& \textbf{NABLA$(k_1,k_2)$ [This paper]}
& \textbf{Yes}
& \textbf{bilateral block EMA}
& $\boldsymbol{\mathcal O(mn+mk_1+nk_2)}$
& \shortstack[l]{
    $\boldsymbol{\widetilde{\mathcal O}(T^{-1/2} +\sigma^{1/2}T^{-1/4})}$
  }
\\ \hline

\end{tabular}%
}

\vspace{3pt}

\begin{minipage}{\TableScale\textwidth}
\footnotesize
\raggedright
``Noise-adaptive'' means that one theorem, without tuning to $\sigma$,
interpolates between the stochastic and zero-noise regimes. A dash means
that the cited work does not establish a directly comparable theorem; it
does not indicate empirical failure. ${}^{\dagger}$The displayed Muon rate
uses a diminishing learning rate and an exponentially increasing batch
size, and controls the best-iterate expected gradient norm rather than the
averaged-gradient criterion used here. For a $d\times d$ factor with block
size $k$, the block statistic has $\lceil d/k\rceil$ principal blocks and
requires $\mathcal O(dk)$ storage; the final block may be smaller then
$k\times d$.
\end{minipage}

\endgroup
\end{table*}


This paper proposes NABLA, a closed-loop adaptive blockwise bilateral preconditioning framework for stochastic optimization. NABLA combines three ingredients: a corrected AdaGrad-Norm coefficient \cite{ward2020adagrad}, an EMA momentum estimator, and matrix-valued left--right preconditioners inspired by
Shampoo-style methods. The corrected coefficient uses both the cumulative squared gradient norm and its running average, so it preserves the stabilizing effect of AdaGrad-Norm while avoiding unnecessarily fast decay. The resulting algorithm is horizon-free and noise-adaptive, and its update is fully determined by the observed trajectory.

Our main contributions are summarized as follows.
\begin{itemize}[leftmargin=14pt,itemsep=3pt,topsep=2pt]
\item We introduce the unified family NABLA$(k_1,k_2)$. For $\X\in\mathbb R^{m\times n}$, the left and right factors contain blocks of sizes $k_1$ and $k_2$, respectively; NABLA$(m,n)$ is the full bilateral realization and NABLA$(1,1)$ is the diagonal realization.

\item We develop a predictable historical-preconditioner analysis. The proof combines monotonicity and ratio control of the closed-loop coefficient, self-normalized momentum bounds, uniform coercivity of the bilateral operator, and an $\mathcal{O}(\beta_t)$ operator-variation estimate.

\item We prove a single nonconvex stationarity bound $\widetilde{\mathcal O}(T^{-1/2}+\sigma^{1/2}T^{-1/4})$ without tuning to $\sigma$. The same theorem yields $\widetilde{\mathcal O}(T^{-1/2})$ when $\sigma=0$, while the block parameterization gives optimizer-state complexity $\mathcal O(mn+mk_1+nk_2)$.

\end{itemize}

The remainder of the paper is organized as follows. Section \ref{sect:related} provides related work. Section \ref{sect:prel} states the assumptions and notation. 
Section \ref{sect:proposed} presents NABLA algorithm. Section \ref{sect:rate} gives the convergence guarantees, Section \ref{sect:experiments} demonstrates the experiment results, and Section \ref{sect:conc} concludes this paper.

\section{Related Work}\label{sect:related}

NABLA is a structured adaptive method for stochastic nonconvex matrix
optimization. We briefly review the most closely related optimizers according
to the geometry retained by their state; Table~\ref{tab:structured-optimizers}
summarizes their preconditioners, memory requirements, and theoretical scope.

Coordinate-wise methods maintain scalar statistics for individual parameters.
AdaGrad accumulates squared gradients and provides data-dependent regret
guarantees in online and stochastic convex optimization~\citep{duchi2011adagrad}.
Adam replaces cumulative statistics by exponential moving averages of the
first and second moments~\citep{kingma2015adam}. Several variants such as AMSGrad \citep{reddi2018convergence} and AdaBelief \citep{zhuang2020adabelief} have been proposed to improve the theoretical convergence behavior of Adam-type methods. OptEMA further makes the EMA
coefficients trajectory-adaptive and proves a single noise-adaptive nonconvex
stationarity guarantee that recovers the deterministic rate as the noise
vanishes~\citep{yuan2026optema}. These methods require only linear optimizer
memory, but their diagonal statistics do not retain correlations between the
row and column modes of a matrix parameter.

One-sided methods retain one mode-sized matrix factor. Muon \cite{kim2026convergence,chang2025convergence} orthogonalizes a
matrix momentum through a polar factor and has been scaled to language-model
training~\citep{liu2025muon}; its recent nonconvex analysis relies on prescribed
learning-rate and batch-size schedules~\citep{nagashima2026muon}. A one-sided
variant of Shampoo has also been analyzed for online regret and offline convex
optimization~\citep{xie2025structured}. OptMuon combines the same one-sided
polar geometry with a closed-loop, noise-adaptive schedule~\citep{yuan2026optmuon}.
For an $m\times n$ matrix with $m\ge n$, this family stores an $n\times n$
factor, but it does not maintain separate factors for both matrix modes.

Two-sided methods retain distinct row and column geometry. Shampoo introduced Kronecker-factored left and right preconditioners with stochastic convex guarantees~\citep{gupta2018shampoo}, and subsequent work improved matrix-root computation and distributed implementation~\citep{anil2020scalable}. Kronecker-factored natural-gradient methods such as K-FAC also exploit structured curvature information for neural network optimization \citep{martens2015optimizing}. KrADagrad avoids inverse matrix roots while preserving bilateral Kronecker structure~\citep{mei2023kradagrad}, whereas SOAP performs an Adam-like second-moment update in the evolving Shampoo eigenbasis~\citep{vyas2025soap}. NABLA$(k_1,k_2)$ differs in two respects: it makes the averaging coefficient of both factors trajectory adaptive, and it exposes block granularity as part of the method. Its structured statistics require $\mathcal O(mk_1+nk_2)$ memory, with the full Shampoo-style geometry recovered at $(k_1,k_2)=(m,n)$.

\section{Preliminaries} \label{sect:prel}

We first fix the notation used in the algorithm and analysis. For matrices,
$\langle\A,\B\rangle:=\operatorname{tr}(\A^\top\B)$ and
$\|\cdot\|_{\fro}$ denote the Frobenius inner product and norm,
respectively, while $\|\cdot\|_{\rm sp}$ denotes the spectral norm.
Let $\mathcal F_t:=\sigma(\xi_1,\ldots,\xi_t)$, where $\sigma(\cdot)$
denotes the sigma-algebra generated by its arguments, and set
$\mathcal F_0:=\{\emptyset,\Omega\}$. We write
$\mathbb E_t[\cdot]:=\mathbb E[\cdot\mid\mathcal F_{t-1}]$.
All unconditional expectations are taken with respect to the randomness of
the stochastic trajectory.

For a symmetric matrix $\A\in\mathbb R^{d\times d}$ and a block size $k\le d$, let $\mathcal P_k(\A)$ retain the diagonal principal blocks induced by consecutive index sets of size $k$ and set all off-block entries to zero. To simplify notation, we assume $k\mid d$, i.e., the dimension $d$ is a multiple of the block size $k$. Accordingly, NABLA$(k_1,k_2)$ uses $m/k_1$ left blocks of size $k_1\times k_1$ and $n/k_2$ right blocks of size $k_2\times k_2$. In particular, NABLA$(m,n)$ denotes the full-matrix configuration, since its block sizes span the entire row and column dimensions of an $m\times n$ matrix.

\begin{assumption}[Lower-Bounded Objective]\label{ass:F}There exists a constant $f_* > -\infty$ such that $f(\X)\ge f_*$ for all $\X\in\mathbb{R}^{m\times n}$.
\end{assumption}

\begin{assumption}[Unbiased Stochastic Gradient with Bounded Variance]\label{ass:sigma} There exists a constant $\sigma\ge0$ such that, for every $\X\in\mathbb R^{m\times n}$, $\mathbb E\!\left[\nabla f(\X;\xi)\right]=\nabla f(\X)$, $\mathbb E \big[\|\nabla f(\X;\xi)-\nabla f(\X)\|_{\fro}^2\big]\le\sigma^2$.
\end{assumption}

\begin{assumption}[Smoothness of the Objective]\label{ass:f} The objective $f$ is $L$-smooth with respect to the Frobenius norm:
\[
\ts \|\nabla f(\X)-\nabla f(\Y)\|_{\fro} \le L\|\X-\Y\|_{\fro},~~\forall\,\X,\Y\in\mathbb{R}^{m\times n}.
\]
\end{assumption}

\begin{assumption}[Spectral Boundedness of Stochastic Gradients]\label{ass:bounded-gradient}
There exists a constant $G>0$ such that
\[
\ts
\|\nabla f(\X;\xi)\|_{\rm sp}\le G,\qquad \forall\,\X\in\mathbb{R}^{m\times n},\ \forall\,\xi.
\]
Consequently, $\|\nabla f(\X;\xi)\|_{\fro}\le {\hat{G}}:=\sqrt{\min\{m,n\}} G$. 
\end{assumption}

\section{The NABLA Algorithm}
\label{sect:proposed}

NABLA separates scalar trajectory adaptation from matrix-valued geometric
preconditioning. The scalar feedback controls the effective memory of the
second-moment estimators, while the block sizes determine how much row- and
column-space correlation is retained. Algorithm~\ref{alg:main} summarizes the
complete update.

\subsection{Closed-Loop Scalar Adaptivity}

For the stochastic-gradient trajectory
$\{\G_i\}_{i\geq1}$, define
\[
\ts
\rho_t
=
\sqrt{
\frac{
1+\frac{1}{t}\sum_{i=1}^{t}\|\G_i\|_{\fro}^{2}
}{
1+\sum_{i=1}^{t}\|\G_i\|_{\fro}^{2}
}}.
\]
The numerator incorporates the empirical average of the observed
gradient energy, whereas the denominator records its cumulative value.
Consequently, $\rho_t$ provides a horizon-free feedback signal determined
entirely by the realized optimization trajectory, without requiring prior
knowledge of the stochastic-gradient noise level.

\begingroup
\def\[{\begin{lrbox}{\optalgbox}\(\displaystyle}
\def\]{\)\end{lrbox}%
\ifdim\wd\optalgbox>.99\columnwidth
  \resizebox{0.99\columnwidth}{!}{\usebox{\optalgbox}}%
\else
  \usebox{\optalgbox}%
\fi}
\begin{figure}[!t]
\refstepcounter{optalg}
\centering
\[
\begin{array}{l}
\hline
\mathbf{Algorithm\ \theoptalg}\ \text{NABLA$(k_1,k_2)$} \\
\hline
01:~\textbf{Input: } \text{base learning rate } \theta>0,\ 
\text{initial point } \X_1\in\mathbb R^{m\times n},\\
~~~~~~~~~ k_1 \in [m],\ k_2 \in [n],\ p,q\ge0\ \text{with}\ p+q=1/2, \varepsilon>0,\ 
\alpha \in(0,1].\\
02:~\text{(Default: } \varepsilon=10^{-6}, \alpha=0.08, \theta=1, \text{and} ~p=q=\tfrac{1}{4}\text{).}\\
03:~\textbf{Initialize: } \X_0=\X_1,~\M_0=\mathbf{0},~\mathbf{L}_0=\mathbf{0},~\R_0=\mathbf{0},~\rho_0=\gamma_0=1.\\
04:~\textbf{for } t=1\,\textbf{to}\,T\,\textbf{do}\\
05:~~~\text{Draw a mini-batch } \xi_t
\text{ and compute}~\G_t:=\nabla f(\X_t;\xi_t) .\\
06:~~~\text{Set}~\beta_t = \gamma_{t-1}.\\
07:~~~\M_t=(1-\alpha)\M_{t-1}+\alpha\G_t.\\
08:~~~\mathbf{L}_t=(1-\beta_t)\mathbf{L}_{t-1}+\beta_t \mathcal{P}_{k_1} (\G_t\G_t^\top).\\
09:~~~\R_t=(1-\beta_t)\R_{t-1}+\beta_t\mathcal{P}_{k_2}(\G_t^\top\G_t).\\
10:~~~\text{Set~}\tilde{\mathbf{L}}_t:=\mathbf{L}_t+\varepsilon\I,~\tilde{\mathbf{R}}_t:=\R_t+\varepsilon\I,~\gamma_t= (\tfrac{
1+\tfrac{1}{t}\sum_{j=1}^{t}\|\M_j\|_{\fro}^{2}}{
1+\sum_{j=1}^{t}\|\M_j\|_{\fro}^{2}})^{1/2}. \\
11:~~~\text{Let~}\X_{t+1} = \X_t - \theta\gamma_t\, \tilde{\mathbf{L}}_t^{-p} \M_t \tilde{\mathbf{R}}_t^{-q} .\\
12:~\textbf{end for}\\
\hline
\end{array}
\]
\label{alg:main}
\end{figure}
\endgroup

NABLA uses the scalar feedback only to regulate the memory of the
matrix-valued second-moment estimators. In particular, the first-moment
coefficient remains fixed at $\alpha$, while the common second-moment
coefficient is chosen as $\beta_t=\rho_{t-1}$. For a generic recursion of the form
\[
\mathbf V_t
=
(1-\beta_t)\mathbf V_{t-1}
+
\beta_t\mathbf S_t,
\]
a larger value of $\beta_t$ assigns more weight to the current
mini-batch statistic $\mathbf S_t$ and therefore produces a more responsive
estimator. Conversely, a smaller value of $\beta_t$ places greater weight on
the historical state and increases its effective memory. The resulting
adaptation changes continuously with the observed gradient trajectory rather
than following a prescribed decay schedule.

The one-step lag in $\beta_t=\rho_{t-1}$ is also important for the analysis.
Before the mini-batch $\xi_t$ is sampled, $\beta_t$ is
$\mathcal F_{t-1}$-measurable. This predictability allows the current
stochastic gradient to be handled through conditional expectation while the
corresponding historical preconditioner remains fixed.

\subsection{Blockwise Bilateral Preconditioning}

NABLA$(k_1,k_2)$ constructs separate second-moment statistics for the row
and column modes of the matrix gradient:
\[
\ts
\begin{aligned}
\mathbf L_t
&=
(1-\beta_t)\mathbf L_{t-1}
+
\beta_t\mathcal P_{k_1}(\G_t\G_t^\top),\\
\R_t
&=
(1-\beta_t)\R_{t-1}
+
\beta_t\mathcal P_{k_2}(\G_t^\top\G_t).
\end{aligned}
\]
The block projections preserve correlations within each selected row or
column block while discarding interactions across different blocks. Under the
divisibility convention adopted above, $\mathbf L_t$ contains $m/k_1$
blocks of size $k_1\times k_1$, and $\R_t$ contains $n/k_2$ blocks of size
$k_2\times k_2$.

The block sizes therefore provide a direct control over the geometric
resolution of the preconditioner. At one extreme, NABLA$(1,1)$ retains only
diagonal row and column statistics. At the other extreme,
NABLA$(m,n)$ retains the complete Gram matrices and corresponds to the full
bilateral realization. Intermediate choices of $(k_1,k_2)$ preserve
within-block matrix geometry while reducing the second-moment storage from
the full-matrix cost to $\mathcal O(mk_1+nk_2)$.

After regularization, define the time-varying bilateral preconditioning
operator
\beq 
\ts \mathcal Q_t(\M):=\widetilde{\mathbf L}_t^{-p}\M\widetilde{\mathbf R}_t^{-q}.
\eeq
The NABLA update can then be written compactly as
\[
\X_{t+1}
=
\X_t-\theta\gamma_t\mathcal Q_t(\M_t).
\]
This representation separates the three roles in the method:
$\M_t$ estimates the gradient direction, $\gamma_t$ controls the overall
update magnitude, and $\mathcal Q_t$ introduces blockwise bilateral
geometry. Changing $(k_1,k_2)$ modifies only the amount of matrix
correlation retained by $\mathcal Q_t$ and its storage cost; it does not
alter the closed-loop scalar feedback or the momentum recursion.

\section{Convergence Guarantees}\label{sect:rate}

In this section, we establish the convergence guarantees of the proposed NABLA  algorithm.

\textbf{Notation.} Throughout this section, we define $\MMM_t:=\sum_{i=1}^{t}\|\M_i\|_{\fro}^{2}$, $\GGG_t:=\sum_{i=1}^{t}\|\G_i\|_{\fro}^{2}$ with $\MMM_0=\GGG_0=0$. Clearly, the variables $\X_t$, $\nabla f(\X_t)$, $\rho_{t-1}$, $\beta_t$, $\gamma_t$, and $\mathcal Q_{t-1}$ are $\mathcal F_{t-1}$-measurable, while $\G_t$, $\M_t$, and $\mathcal Q_t$ are not $\mathcal F_{t-1}$-measurable.

First, we establish some preliminary bounds.
These results characterize the stability of the adaptive stepsize and
bilateral preconditioner, which will be used to derive the descent
inequality.

\begin{lemma} \label{lemma:preliminary}
(Proof in Appendix \ref{app:lemma:preliminary}, Preliminary Bounds) Define ${\hat{G}}:=\sqrt{\min\{m,n\}}G$, $h:=\varpi^{p+q}$, and $C_a:=G^2\varepsilon^{-3/2}$. For every $1\le t\le T$, the following statements hold:

\begin{enumerate}[label=(\alph*)]
\item
$\MMM_T\le \GGG_T$, and $\|\M_t\|_{\fro}\le {\hat{G}}$.

\item $1\le \frac{\gamma_t}{\gamma_{t+1}}
\le \sqrt{\kappa}$, where $\kappa:=1+{\hat{G}}^2$.

\item $\|\mathbf L_t\|_{\rm sp}\le G^2$, $\|\mathbf R_t\|_{\rm sp}\le G^2$.

\item $\varepsilon\I\preceq\widetilde{\mathbf L}_t,\widetilde{\mathbf R}_t
\preceq\varpi\I$, where $\varpi:=\varepsilon+G^2$.

\item $\langle \M_t,\mathcal Q_t(\M_t)\rangle \ge \frac{1}{h}\|\M_t\|_{\fro}^{2}$, $\|\mathcal Q_t(\M_t)\|_{\fro}^2\le \tfrac{1}{\varepsilon}\|\M_t\|_{\fro}^2$.

\item $\|\mathcal Q_t(\M_t)-\mathcal Q_{t-1}(\M_t)\|_{\fro} \le C_a\beta_t\|\M_t\|_{\fro}$.

\end{enumerate}
\end{lemma}

The above bounds ensure the uniform boundedness and variation control of
the preconditioner. We next establish the one-step descent inequality by
decomposing the update into a predictable historical component and a
preconditioner variation error.

\begin{lemma} \label{lemma:one:step:descent}

(Proof in Appendix \ref{app:lemma:one:step:descent},  One-Step Descent Inequality) Let $c_0:=\frac{\theta {\hat{G}}^2}{\sqrt\varepsilon}$, $c_1:=\frac{L\theta^2}{\varepsilon}$, $c_2:=\frac{\theta h C_a^2}{\alpha}$, $c_3:=\frac{L\theta^2}{2\varepsilon}$, $c_4:=\frac{3\theta}{4h}$. For every $t\ge1$,
\beq
&& \ts \mathbb E[f(\X_{t+1})-f(\X_t)]  \le  \mathbb E[ c_3\gamma_t^2\|\M_t\|_{\fro}^{2}] \nn \\
&& \ts + \mathbb E[ \sum_{j=1}^{t}(1-\alpha)^{t-j}(E_j-D_j)] \nn\\
&& \ts + \mathbb E[  (1-\alpha)c_0 \sum_{j=2}^{t}(1-\alpha)^{t-j}(\gamma_{j-1}-\gamma_j)].\nn
\eeq
Here, $E_t:=c_1\gamma_{t-1}^2\|\M_{t-1}\|_{\fro}^{2} + c_2\beta_t^2\gamma_t\|\M_t\|_{\fro}^{2}$, and $D_t:=c_4\alpha\gamma_t\|\nabla f(\X_t)\|_{\fro}^{2}$.
    
\end{lemma}

The one-step descent inequality contains error terms induced by momentum
propagation and preconditioner variation. We next show that these terms
are controlled by logarithmic self-normalized bounds.

\begin{lemma} \label{lemma:log:bound}

(Proof in Appendix \ref{app:lemma:log:bound}, Logarithmical bounds) Let $C_b:=(1+\widehat G^2) \ln (e+\hG^2)$, $C_e:=\frac{c_1\kappa}{\alpha}
+\frac{c_2\kappa^3}{\alpha} + c_3\kappa$, and $C_h:=\frac{ f(\X_1) - f_* + C_e+c_3 C_b+(1-\alpha)c_0/\alpha}{
c_4\alpha}$. Define $H_T:=\sum_{t=1}^{T}\gamma_t\|\nabla f(\X_t)\|_{\fro}^{2}$. Then

(a) $\sum_{t=1}^{T}\gamma_t^2\|\M_t\|_{\fro}^2 \le C_b \ln(e+T)$ almost surely.

(b) $\sum_{t=1}^{T}\sum_{j=1}^{t}(1-\alpha)^{t-j}E_j
\le C_e\ln(1+T\widehat G^2)$.

(c) $\mathbb E[H_T]\le C_h\ln(e+T)$.

\end{lemma}

The above result provides a logarithmic bound on the weighted gradient
energy. The following lemma converts this weighted quantity into the
standard cumulative gradient measure.

\begin{lemma} \label{lemma:Weighted:energy:bridge} (Proof in Appendix \ref{app:lemma:Weighted:energy:bridge}) Let $A_T:=\sum_{t=1}^{T}\|\nabla f(\X_t)\|_{\fro}^{2}$, $Z_T:=\mathbb E\sqrt{1+\GGG_T}$. Let $C_Z=\sqrt{C_h}$. We have: $\mathbb E\left[\sqrt{A_T}\right] \le C_Z Z_T^{1/2}\ln^{1/2}(e+T)$. 
\end{lemma}

Combining the weighted-energy bridge with the stochastic-gradient
decomposition yields the final noise-adaptive convergence guarantee.

\begin{theorem} \label{the:1}
(Noise-adaptive convergence and zero-noise \mbox{optimality}, Proof in Appendix \ref{app:the:1}) Under Assumptions~\ref{ass:F}--\ref{ass:bounded-gradient},
for every $T\ge1$, Algorithm~\ref{alg:main} satisfies
\beq
&& \ts \mathbb E \left[
\frac1T\sum_{t=1}^{T}\|\nabla f(\X_t)\|_{\fro}
\right] \nn\\
&\le& \mathcal O\Bigl(
T^{-1/2}\ln(e+T)  + \sigma^{1/2}T^{-1/4}\ln^{1/2}(e+T)
\Bigr).\nn
\eeq
If $R\sim\operatorname{Unif}\{1,\ldots,T\}$
is independent of the algorithmic randomness and
$\X_{\rm out}:=\X_R$,
then $\mathbb E
\|\nabla f(\X_{\rm out})\|_{\fro} = \mathbb E
\left[
\frac1T
\sum_{t=1}^{T}
\|\nabla f(\X_t)\|_{\fro}
\right]$ and the same bound holds.
\end{theorem}

Theorem~\ref{the:1} establishes a unified convergence guarantee that
interpolates between stochastic and deterministic optimization regimes.
Unlike methods requiring prior knowledge of the noise level, NABLA
automatically adapts to the observed stochastic-gradient trajectory and
achieves the rate $\widetilde{\mathcal O}(T^{-1/2}+\sigma^{1/2}T^{-1/4})$ without tuning to $\sigma$. In particular, when the stochastic gradient noise vanishes,
the bound recovers the deterministic-order rate
$\widetilde{\mathcal O}(T^{-1/2})$, demonstrating that the proposed
closed-loop adaptation does not introduce an additional stochastic
penalty in the zero-noise regime. 

\section{Experiments}\label{sect:experiments}

We evaluated NABLA on two complementary learning problems: image
classification with ResNet-18 on CIFAR-100 and language-model pretraining with
LLaMA-60M on C4. The experiments addressed two questions. First, how does the
optimization quality of NABLA compare with commonly used adaptive and
matrix-preconditioned optimizers? Second, how does block granularity trade
retained matrix geometry for optimizer-state memory and wall-clock cost?
Accordingly, we compared the shape-matched configuration NABLA$(m,n)$ with
NABLA$(64,128)$ on CIFAR-100 and NABLA$(1,1)$ on C4.

\subsection{Experimental Setup}\label{sect:exp-setup}

\paragraph{CIFAR-100 image classification.}
We trained a CIFAR-style ResNet-18~\citep{he2016resnet} on
CIFAR-100~\citep{krizhevsky2009cifar} for 100 epochs. The network uses a
$3\times3$ stride-one stem without max pooling, followed by four residual
stages with $64$, $128$, $256$, and $512$ channels and an
adaptive-average-pooling classifier. Training uses a four-pixel padded random
crop, random horizontal flipping,
TrivialAugmentWide~\citep{muller2021trivialaugment}, and cross entropy with
label smoothing $0.1$. The training and test batch sizes are $128$ and $256$,
respectively. Model parameters and optimizer states are stored in FP32, while
forward and evaluation passes use BF16 autocast without gradient scaling.

All methods use a 600-update linear warm-up. NABLA holds its base scale through
epoch $60$ and subsequently applies cosine decay to a final multiplier of
$0.25$; the baselines use their tuned optimizer-specific schedules. For NABLA,
we set $\theta=7$, $\tau=1$, $\alpha=0.08$, $\varepsilon=10^{-6}$, and
$p=q=1/4$. NABLA$(64,128)$ uses left and right block sizes $k_1=64$ and
$k_2=128$, whereas NABLA$(m,n)$ retains one shape-matched block along each
matrix mode.

\paragraph{LLaMA-60M language modeling.}
We pretrained a 58.1M-parameter LLaMA-style causal language
model~\citep{touvron2023llama} on the English C4
corpus~\citep{raffel2020t5}. The model has hidden size $512$, eight layers,
eight attention heads, eight key--value heads, intermediate size $1376$, and a
$32{,}000$-token vocabulary. Text is tokenized with \texttt{t5-base} and
truncated or padded to length $256$. Each run uses two GPUs, per-GPU batch size
$128$, gradient accumulation $2$, global batch size $512$, FP32 training,
gradient clipping at $1.0$, and a 2,000-update warm-up. Training progress is
reported in processed tokens, and validation is performed at a common schedule
for all compared configurations.

For NABLA, the selected base learning rate is $\theta=0.05$, while the AdamW fallback learning
rate is $0.001$, $\alpha=0.05$, $\varepsilon=10^{-5}$, and $p=q=1/4$.
NABLA$(m,n)$ and NABLA$(1,1)$ use the same momentum recursion, scalar feedback,
inverse-factor postprocessing, and update stabilization. Their only structural
difference is block granularity: NABLA$(m,n)$ retains complete left and right
Gram statistics, whereas NABLA$(1,1)$ retains their diagonal entries.

\paragraph{Compared methods and repeated runs.}
On CIFAR-100, we compared NABLA$(m,n)$ and NABLA$(64,128)$ with
AdamW~\citep{loshchilov2019adamw},
Shampoo~\citep{gupta2018shampoo}, and
Muon~\citep{liu2025muon}. AdamW uses learning rate $10^{-3}$ and weight decay
$5\times10^{-4}$. Shampoo uses matrix learning rate $0.5$, momentum $0.9$,
regularizer $10^{-6}$, and refreshes its inverse roots every 10 updates. Muon
uses learning rate $10^{-2}$ and momentum $0.95$.

The CIFAR-100 comparison uses six paired random seeds. For each seed, all
optimizers share the model initialization, data order, augmentation
randomness, training horizon, and evaluation epochs. The C4 comparison uses
three paired seeds for NABLA$(m,n)$ and NABLA$(1,1)$. The Shampoo, Muon, and
AdamW trajectories are included as scale-matched optimizer references. Within
each benchmark, wall-clock measurements use the same hardware and timing
definition. All NABLA runs implement Algorithm~\ref{alg:main}; benchmark-level
warm-up and learning-rate schedules modify only the base learning-rate
multiplier and do not alter its feedback or block-preconditioner recursions.
Complete implementation settings are provided in
Appendix~\ref{supp:sec:implementation}.
\paragraph{Sensitivity-study protocols.}
We complemented the paired-seed comparisons with two diagnostic sensitivity
studies. On CIFAR-100, we performed a fixed-seed sweep of the NABLA$(64,128)$
base scale over
\(\theta\in[10^{-1},2\times10^{1}]\) using the same ResNet-18 data and
augmentation pipeline. For each setting, we recorded the minimum training and
test losses and the maximum training and test accuracies attained within the
fixed sweep horizon. This sweep was used to identify a stable operating region
and was not pooled with the six-seed optimizer comparison.

On C4, we performed a strict base-learning-rate sweep over
\(\theta\in[10^{-4},10^{-1}]\), holding the remaining training settings fixed,
and recorded the evaluation loss at \(10{,}000\) updates for NABLA$(m,n)$ and
NABLA$(1,1)$. At the selected value \(\theta=0.05\), we additionally measured
validation perplexity at \(8{,}000\), \(9{,}000\), and \(10{,}000\) updates
for the same two block granularities. The first experiment evaluates
learning-rate sensitivity, whereas the second isolates the effect of replacing
the complete bilateral factors with diagonal statistics.
\subsection{Training Dynamics and Wall-Clock Cost}
\label{sect:exp-dynamics}

\begin{figure*}[t!]
\centering
\includegraphics[width=\textwidth]
{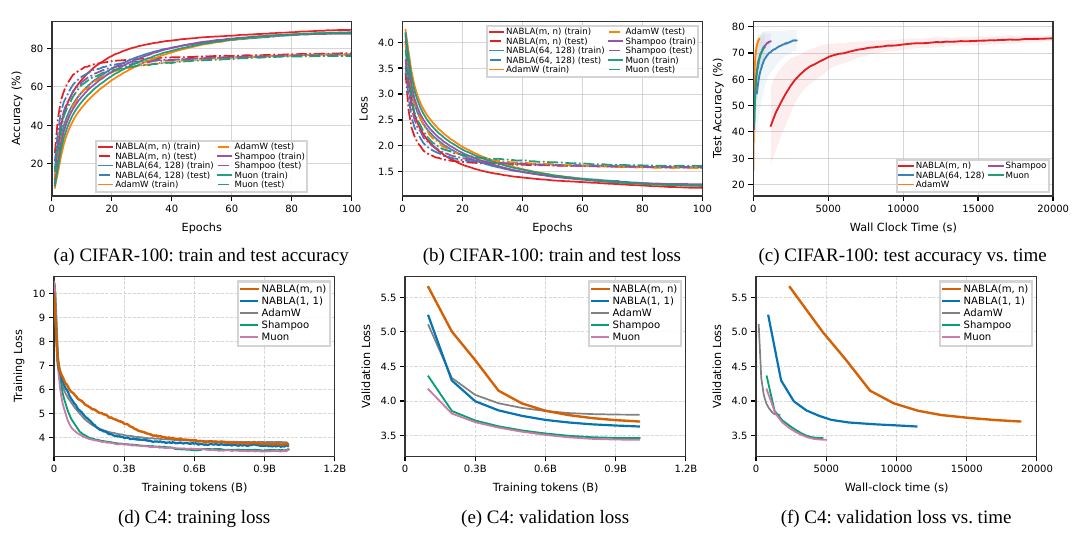}
\caption{Training dynamics of NABLA configurations and baseline optimizers.
The top row reports ResNet-18 results on CIFAR-100:
\textbf{(a)} training and test accuracy versus epoch,
\textbf{(b)} training and test loss versus epoch, and
\textbf{(c)} test accuracy versus wall-clock time.
The bottom row reports LLaMA-60M pretraining on C4:
\textbf{(d)} training loss versus processed tokens,
\textbf{(e)} validation loss versus processed tokens, and
\textbf{(f)} validation loss versus wall-clock time.
The CIFAR-100 comparison evaluates NABLA$(m,n)$ and NABLA$(64,128)$, whereas
the C4 comparison evaluates NABLA$(m,n)$ and NABLA$(1,1)$. AdamW, Shampoo,
and Muon are included as baseline optimizers in both benchmarks.}
\label{fig:training-dynamics}
\end{figure*}

\paragraph{CIFAR-100 optimization.}
Figures~\ref{fig:training-dynamics}(a) and
\ref{fig:training-dynamics}(b) show that all five optimizers reduce the
classification loss rapidly during the early epochs. NABLA$(m,n)$ reaches the
highest training-accuracy trajectory and the lowest training-loss trajectory,
while NABLA$(64,128)$ follows the same overall convergence pattern with
substantially smaller preconditioner blocks. The test curves are more tightly
clustered than the training curves, indicating that the principal empirical
difference between the two NABLA configurations lies in the geometry and cost
of their optimizer states rather than in a large separation in predictive
accuracy.

The wall-clock view in Figure~\ref{fig:training-dynamics}(c) changes the
interpretation of the epoch-wise comparison. The shape-matched
NABLA$(m,n)$ configuration requires substantially more elapsed time because it
maintains and transforms larger matrix factors. NABLA$(64,128)$ reaches its
terminal accuracy region much earlier and has a wall-clock profile closer to
the baseline optimizers. Thus, blockwise preconditioning retains most of the
epoch-wise optimization behavior while reducing the computational cost of the
shape-matched realization.

\paragraph{C4 language-model pretraining.}
Figures~\ref{fig:training-dynamics}(d) and
\ref{fig:training-dynamics}(e) show that all five optimizers steadily reduce
the training and validation losses as the number of processed tokens
increases. Within the NABLA family, NABLA$(1,1)$ consistently follows a lower
loss trajectory than NABLA$(m,n)$ after the early training stage. At the final
reported checkpoint, the validation losses are approximately \(3.62\) for
NABLA$(1,1)$ and \(3.70\) for NABLA$(m,n)$. NABLA$(1,1)$ also improves upon
the displayed AdamW endpoint, whereas Shampoo and Muon attain lower final
validation losses in the present implementation. The experiment therefore
does not indicate that retaining complete left and right Gram statistics
necessarily improves optimization at this model scale; the diagonal
realization provides the stronger result among the two tested NABLA
configurations.

The wall-clock comparison in Figure~\ref{fig:training-dynamics}(f) reinforces
this distinction. The displayed NABLA$(1,1)$ trajectory reaches its terminal
loss region in approximately \(11{,}500\) seconds, compared with approximately
\(19{,}000\) seconds for NABLA$(m,n)$, corresponding to roughly \(40\%\) less
elapsed time in the present implementation. Moreover, NABLA$(1,1)$ reaches a
lower terminal validation loss despite using the less expensive diagonal
statistics. Shampoo and Muon remain faster and reach lower validation losses,
so the result should not be interpreted as universal wall-clock dominance.
Instead, it shows that block granularity materially affects the quality--cost
trade-off within NABLA and that the diagonal realization is preferable for
the tested LLaMA-60M/C4 setting.

\subsection{Learning-Rate and Block-Granularity Sensitivity}
\label{sect:exp-sensitivity}

\paragraph{CIFAR-100 base-scale sensitivity.}
Figure~\ref{fig:cifar-lr-sensitivity} shows that the selected base scale lies
inside a broad stable region rather than at an isolated optimum. At
\(\theta=0.1\), NABLA$(64,128)$ attains a minimum training loss of approximately
\(1.70\), a minimum test loss of approximately \(1.82\), and maximum training
and test accuracies of approximately \(74.6\%\) and \(69.1\%\), respectively.
Performance improves rapidly as the base scale increases. For
\(\theta\gtrsim3\), the minimum training loss remains approximately within
\(1.36\)--\(1.38\), while the maximum training accuracy remains within
approximately \(84.6\%\)--\(85.0\%\). Over the same region, the minimum test
loss varies within approximately \(1.60\)--\(1.64\), and the maximum test
accuracy remains around \(75\%\)--\(76\%\). The main-experiment choice
\(\theta=7\) therefore falls within the observed plateau. The fluctuations at
larger scales also show that the sweep does not support a unique globally
optimal learning rate; rather, it identifies a range in which the tested
configuration is not sensitive to fine-grained tuning.

\begin{figure*}[t!]
\centering
\includegraphics[width=0.70\textwidth]
{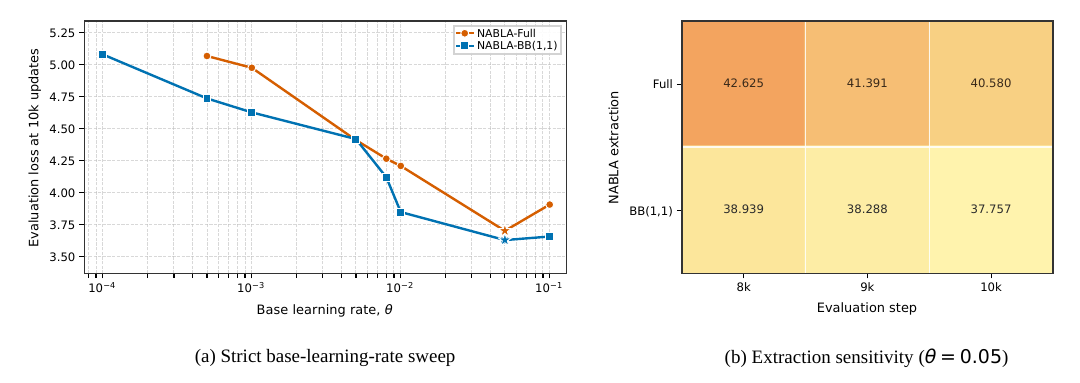}
\caption{Base-learning-rate and block-granularity sensitivity of NABLA on C4
with LLaMA-60M. \textbf{(a)} Evaluation loss at \(10{,}000\) updates for
NABLA$(m,n)$ and NABLA$(1,1)$ across the evaluated base-learning-rate range.
\textbf{(b)} Validation perplexity of NABLA$(m,n)$ and NABLA$(1,1)$ at
\(8{,}000\), \(9{,}000\), and \(10{,}000\) updates under the selected base
learning rate \(\theta=0.05\). Lower values indicate better performance in
both panels.}
\label{fig:c4-lr-extraction-sensitivity}
\end{figure*}

\paragraph{C4 learning-rate and block-granularity sensitivity.}
Figure~\ref{fig:c4-lr-extraction-sensitivity}(a) shows that the two NABLA
configurations respond differently to the base learning rate. Increasing
\(\theta\) from the smallest tested values substantially reduces the
evaluation loss at \(10{,}000\) updates for both configurations. The lowest
displayed losses are attained near \(\theta=0.05\): NABLA$(m,n)$ reaches
approximately \(3.71\), whereas NABLA$(1,1)$ reaches approximately \(3.63\).
Increasing the base learning rate further to \(\theta=0.1\) raises the loss
of NABLA$(m,n)$ to approximately \(3.91\), while NABLA$(1,1)$ remains near
\(3.66\). Thus, \(\theta=0.05\) lies near the best-performing region of the tested sweep, and the diagonal realization is less sensitive to the largest evaluated learning rate.

Figure~\ref{fig:c4-lr-extraction-sensitivity}(b) confirms the advantage of
NABLA$(1,1)$ at the selected base learning rate. From \(8{,}000\) to
\(10{,}000\) updates, the validation perplexity of NABLA$(m,n)$ decreases
from \(42.625\) to \(40.580\), whereas that of NABLA$(1,1)$ decreases from
\(38.939\) to \(37.757\). The corresponding perplexity gaps are \(3.686\),
\(3.103\), and \(2.823\) at \(8{,}000\), \(9{,}000\), and \(10{,}000\)
updates, respectively. Although the gap narrows as training progresses, it
remains visible at every reported checkpoint. Under this model scale and
training configuration, diagonal block granularity therefore improves both
the optimization endpoint and the associated computational cost relative to
the shape-matched realization.

\subsection{Block Granularity and Optimizer-State Memory}
\label{sect:exp-memory}

To isolate persistent optimizer-state storage, we profiled the two CIFAR-100
NABLA configurations after lazy state initialization using the same ResNet-18
and one synthetic $128\times3\times32\times32$ training batch. Profiling used
FP32 parameters, BF16 autocast, PyTorch 2.9.1 with CUDA 12.8, and an NVIDIA RTX
PRO 6000 Blackwell Server Edition.

\begin{table}[t!]
\centering
\caption{Measured memory of the two NABLA configurations on the CIFAR-100
ResNet-18 setup. ``State'' denotes persistent optimizer-state tensors, whereas
``Peak'' includes the complete training-step allocation. Lower is better.}
\label{tab:nabla-memory}
\setlength{\tabcolsep}{4.5pt}
\renewcommand{\arraystretch}{1.10}
\begin{tabular}{lrr}
\toprule
Configuration &
\shortstack{State\\(MiB) $\downarrow$} &
\shortstack{Peak\\(MiB) $\downarrow$}\\
\midrule
NABLA$(m,n)$      & 760.014 & 1621.962\\
NABLA$(64,128)$   & \textbf{76.282} & \textbf{391.649}\\
\midrule
Reduction         & $9.96\times$ & $4.14\times$\\
\bottomrule
\end{tabular}
\end{table}

As shown in Table~\ref{tab:nabla-memory}, replacing the shape-matched factors
with $64\times128$ blocks reduces persistent optimizer-state storage from
$760.014$ MiB to $76.282$ MiB, a $9.96\times$ reduction. Peak allocated memory
decreases from $1621.962$ MiB to $391.649$ MiB, corresponding to a
$4.14\times$ reduction. The smaller peak-memory ratio is expected because peak
allocation also contains model parameters, activations, gradients, input
tensors, and temporary matrix-computation buffers. We therefore use
$9.96\times$ only as the persistent-state reduction and do not present it as a
total GPU-memory reduction.

\subsection{Evidence Boundaries}\label{sect:exp-efficiency}

Five considerations delimit the empirical conclusions. First, the CIFAR-100
optimizers use tuned optimizer-specific learning-rate schedules; their
epoch-wise comparison is therefore a comparison between practical training
configurations rather than an ablation under an identical schedule. Second,
the controlled C4 conclusion concerns the paired
NABLA$(m,n)$--NABLA$(1,1)$ comparison. Shampoo, Muon, and AdamW provide
scale-matched context, but they should not be interpreted as paired statistical
comparisons unless they are rerun under the same seed protocol.

Third, the CIFAR-100 learning-rate sweep in
Figure~\ref{fig:cifar-lr-sensitivity} is a fixed-seed diagnostic and reports
within-run extrema. It supports identification of a broad stable operating
region but does not quantify between-seed uncertainty. Fourth, the C4
learning-rate sweep reports the evaluation loss at \(10{,}000\) updates over
the tested range of base learning rates, while the accompanying heatmap
compares only the two endpoint granularities, NABLA$(m,n)$ and NABLA$(1,1)$,
at \(\theta=0.05\). These results support the local conclusion that the
diagonal realization is more favorable under the tested configuration, but
they do not establish a monotone relationship between block size and
optimization quality across other model scales, learning rates, or training
horizons.

Finally, wall-clock results are hardware- and implementation-dependent. They
establish the cost ordering observed in the present implementation but do not
imply universal runtime dominance. Evaluation on additional architectures,
model scales, hardware platforms, learning-rate neighborhoods, and intermediate
block sizes remains necessary to characterize broader scaling behavior.

\section{Conclusion}
\label{sect:conc}

We introduced NABLA$(k_1,k_2)$, a unified family of closed-loop bilateral block-preconditioned methods. Its fixed first-moment EMA is coupled to trajectory-adaptive left and right second-moment statistics, while $(k_1,k_2)$ determines the retained matrix geometry and state memory. The analysis proceeds through predictable coefficient control, self-normalized momentum bounds, uniform bilateral coercivity, operator variation, historical descent, and a final noise closure. The resulting theorem is horizon-free and noise-adaptive, and it recovers the deterministic-order $\widetilde{\mathcal O}(T^{-1/2})$ gradient-norm rate when the stochastic-gradient noise vanishes.

\clearpage
\normalem

\bibliography{mybib}
\appendix

\clearpage
\onecolumn
{\Huge Appendix}
\appendix

{ The appendix is organized as follows.

Appendix~\ref{supp:sec:support} collects supporting lemmas. 

Appendix~\ref{supp:sect:rate} presents the proofs from Section \ref{sect:rate}.

Appendix~\ref{supp:sec:implementation} provides additional implementation details.}

\section{Notations and Supporting Lemmas}\label{supp:sec:support}

\subsection{Notations}

Throughout the paper, we use the following notation. Bold lowercase letters denote vectors, while uppercase letters denote real-valued matrices.

\begin{itemize}[leftmargin=12pt,itemsep=0.2ex]

\item $[n]$: The set $\{1,2,...,n\}$.

\item $\|\mathbf{x}\|$: Euclidean norm, defined as $\|\mathbf{x}\|=\|\mathbf{x}\|_2 = \sqrt{\la \mathbf{x},\mathbf{x}\ra}$.

\item $\la \mathbf{a},\mathbf{b}\ra$ : Standard inner product, given by $\la \mathbf{a},\mathbf{b}\ra =\sum_{i}{\mathbf{a}_{i}\mathbf{b}_{i}}$.

\item $\mathbf{a} \leq \alpha$: For $\mathbf{a}\in\Rn^n$ and scalar $\alpha\in\Rn$, this means $\mathbf{a}_i\leq \alpha$ for all $i\in[n]$.

\item $\EEE[X]$: Expectation of the random variable $X$.

\item $\text{Var}(X)$: Variance of the random variable $X$.

\item $\{A_i\}_{i=1}^{\infty}$, $\{B_i\}_{i=1}^{\infty}$: sequences indexed by non-negative integers.

\item $\A^\top$ : the transpose of the matrix $\A$.

\item $\|\A\|_{\fro}$: Frobenius norm of the matrix $\A$.

\item $\|\A\|_{{\rm sp}}$: Spectral norm of the matrix $\A$.

\item $\X\succeq \zero$: The matrix $\X$ is a symmetric and positive semi-definite matrix.

\item $k\mid d$: The dimension $d$ is a multiple of the block size $k$.

\end{itemize}

\subsection{Supporting Lemmas}
We offer a set of helpful lemmas, each of which stands independently of context and specific methodology

\begin{lemma}[Self-Normalized Summation]\label{lemma:selfnorm} Let $b_1,\ldots,b_T\ge0$ and $B_t:=\sum_{i=1}^{t}b_i$. We have the following results.

\begin{enumerate}[label=(\alph*)]

\item It holds that $\sum_{t=1}^{T}\frac{b_t}{1+B_t}\le\ln(1+B_T)$.

\item If additionally $b_t\le B$ for every $t$, then $\sum_{t=1}^{T}\frac{b_t}{1+B_{t-1}} \le(1+B)\ln(1+B_T)$.

\end{enumerate}

\begin{proof}
\noindent\textbf{Part (a).} We derive:
\beq
\ts \sum_{t=1}^{T}\frac{b_t}{1+B_t} ~~\le~~ \sum_{t=1}^{T} \left(\ln(1+B_t)-\ln(1+B_{t-1})\right) ~~=~~ \ln(1+B_T).\nn
\eeq

\noindent\textbf{Part (b).} We derive:
\beq
\ts \sum_{t=1}^{T}\frac{b_t}{1+B_{t-1}} ~~\xxxx{1}{\le}~~ (1+B)\sum_{t=1}^{T}\frac{b_t}{1+B_{t}}~~\leq~~(1+B)\ln(1+B_T),\nn
\eeq
\noi where step \step{1} uses $1+B_t\le(1+B)(1+B_{t-1})$.

\end{proof}

\end{lemma}

\begin{lemma} \label{lemma:closed-loop-coefficient} (Closed-loop Coefficient Bound) Let $\GGG_t:=\sum_{i=1}^{t}\|\G_i\|_{\fro}^{2}$, $\GGG_0:=0$, and suppose that $\|\G_i\|_{\fro}\le {\hat{G}}$ for every $i$.
Set $\beta_1=1$, and, for every $t\ge2$, let $\beta_{t-1}
:= \big(\frac{1+\GGG_{t}/t}
     {1+\GGG_{t}}\big)^{1/2}$. Then, for every $t\ge1$,
\beq
\ts \beta_t^2 \le \frac{\kappa}{1+\GGG_{t-1}},~~~~\text{where}~~~~\kappa:=1+{\hat{G}}^2.
\eeq

\begin{proof}
For $t=1$, the result follows from $\beta_1=1$, $\GGG_0=0$, and
$\kappa\ge1$. For $t\ge2$, Assumption~\ref{ass:bounded-gradient} gives
\[
\frac{\GGG_{t-1}}{t-1}
=
\frac{1}{t-1}\sum_{i=1}^{t-1}\|\G_i\|_{\fro}^{2}
\le {\hat{G}}^2.
\]
Therefore,
\[
\beta_t^2
=
\frac{1+\GGG_{t-1}/(t-1)}
     {1+\GGG_{t-1}}
\le
\frac{1+{\hat{G}}^2}
     {1+\GGG_{t-1}}
=
\frac{\kappa}
     {1+\GGG_{t-1}}.
\]

\end{proof}

\end{lemma}

\begin{lemma}[Inverse Fractional-Power Perturbation]\label{lemma:inverse-power}
Let $\A,\B\succeq\varepsilon\I$ and $s\in[0,1]$. Then
\[
\|\A^{-s}-\B^{-s}\|_{\rm sp}
\le s\varepsilon^{-s-1}\|\A-\B\|_{\rm sp}.
\]

\begin{proof}
The case $s=0$ is immediate. For $s\in(0,1)$, use
\[
\A^{-s}=\frac{\sin(\pi s)}{\pi}\int_0^\infty
\nu^{-s}(\A+\nu\I)^{-1}\,d\nu
\]
and the resolvent identity. Since
$\|(\A+\nu\I)^{-1}\|_{\rm sp},\|(\B+\nu\I)^{-1}\|_{\rm sp}\le(\varepsilon+\nu)^{-1}$,
\[
\|\A^{-s}-\B^{-s}\|_{\rm sp}
\le
\frac{\sin(\pi s)}{\pi}\|\A-\B\|_{\rm sp}
\int_0^\infty\frac{\nu^{-s}}{(\varepsilon+\nu)^2}\,d\nu.
\]
The integral equals $s\pi\varepsilon^{-s-1}/\sin(\pi s)$.
\end{proof}
\end{lemma}

\begin{lemma}[Fixed-EMA denominator coupling]\label{lemma:ema-denominator}
Let $\M_t=(1-\alpha)\M_{t-1}+\alpha\G_t$, $\M_0=0$, and $\GGG_t=\sum_{i=1}^{t}\|\G_i\|_{\fro}^{2}$. Then, 
\beq
\sum_{t=1}^{T}\frac{\|\M_t\|_{\fro}^{2}}{1+\GGG_{t-1}}
\le\kappa\ln(1+\GGG_T),~~~~\text{where}~~~~\kappa=1+{\hat{G}}^2.
\eeq

\begin{proof} Define $\GGG_t=\sum_{i=1}^{t}\|\G_i\|_{\fro}^{2}$.

\noi The momentum recursion gives
\[
\M_t=\alpha\sum_{i=1}^{t}(1-\alpha)^{t-i}\G_i,
\]
and convexity of the squared Frobenius norm yields
\[
\|\M_t\|_{\fro}^{2}
\le
\alpha\sum_{i=1}^{t}(1-\alpha)^{t-i}\|\G_i\|_{\fro}^{2}.
\]
Therefore,
\beq
\ts \sum_{t=1}^{T}\frac{\|\M_t\|_{\fro}^{2}}{1+\GGG_{t-1}}
&\le& \ts \alpha\sum_{i=1}^{T}\|\G_i\|_{\fro}^{2}
\sum_{t=i}^{T}\frac{(1-\alpha)^{t-i}}{1+\GGG_{t-1}}~~\le~~\sum_{i=1}^{T}\frac{\|\G_i\|_{\fro}^{2}}{1+\GGG_{i-1}} \nn \\
&\xxxx{1}{\leq}& \ts \kappa\sum_{i=1}^{T}\frac{\|\G_i\|_{\fro}^{2}}{1+\GGG_i}
~~\xxxx{2}{\le}~~\kappa\ln(1+\GGG_T),\nn 
\eeq
where step \step{1} uses $\|\G_i\|_{\fro}^{2}\le {\hat{G}}^2$; step \step{2} uses Lemma~\ref{lemma:selfnorm}.
\end{proof}
\end{lemma}

\section{Proof for Section \ref{sect:rate}}\label{supp:sect:rate}


\subsection{Proof of Lemma \ref{lemma:preliminary}}
\label{app:lemma:preliminary}

\begin{proof} Define $\kappa:=1+\hG^2$.

\paragraph{Part (a-i).}

The momentum recursion gives
\beq
\ts \M_t = \alpha\sum_{i=1}^{t} (1-\alpha)^{t-i}\G_i.\nn
\eeq

Because the coefficients form a sub-probability distribution, $\sum_{i=1}^{t}
\alpha(1-\alpha)^{t-i} = 1-(1-\alpha)^t
\le1$, convexity of the squared Frobenius norm yields $\|\M_t\|_{\fro}^{2}
\le \alpha \sum_{i=1}^{t}
(1-\alpha)^{t-i}
\|\G_i\|_{\fro}^{2}$. 

Therefore, we derive:
\beq
\MMM_T &=& \ts \sum_{t=1}^{T}\|\M_t\|_{\fro}^{2} ~~\le~~ \alpha
\sum_{t=1}^{T} \sum_{i=1}^{t} (1-\alpha)^{t-i} \|\G_i\|_{\fro}^{2} \nn \\
&=& \ts \alpha \sum_{i=1}^{T} \|\G_i\|_{\fro}^{2} \sum_{t=i}^{T}(1-\alpha)^{t-i} ~~\le~~ \ts \sum_{i=1}^{T}\|\G_i\|_{\fro}^{2} = \GGG_T. \nn
\eeq

\paragraph{Part (a-ii).} Assumption~\ref{ass:bounded-gradient} implies $\|\G_t\|_{\fro} \le \sqrt{\min\{m,n\}}\,G = {\hat{G}}$. Hence,
\beq
\|\M_t\|_{\fro}^{2} ~~\le~~ \alpha {\hat{G}}^2 \sum_{i=1}^{t}(1-\alpha)^{t-i}
~~=~~ (1-(1-\alpha)^t){\hat{G}}^2
\le {\hat{G}}^2. \nn
\eeq
\noi This leads to:
\beq
\|\M_t\|_{\fro}\le {\hat{G}}. \nn
\eeq

\paragraph{Part (b).} Since $\gamma_t=(1+\MMM_{t-1})^{-1/2}$, the sequence $\{\gamma_t\}$ is non-increasing. Therefore, $\frac{\gamma_t}{\gamma_{t+1}}\ge1$. Additionally,
\beq
\ts \frac{\gamma_t}{\gamma_{t+1}} ~~=~~
\sqrt{\frac{1+\MMM_t}{1+\MMM_{t-1}}}
~~=~~ \sqrt{1+
\frac{\|\M_t\|_{\fro}^2}
{1+\MMM_{t-1}}}
~~\le~~ \sqrt{1+\|\M_t\|_{\fro}^{2}} ~~\leq ~~\sqrt{1 + \hG^2}~~:=~~\sqrt{\kappa}.\nn
\eeq

\paragraph{Part (c).} Define
\[
S_t^L=\mathcal P_{k_1}(\G_t\G_t^\top),
\qquad
S_t^R=\mathcal P_{k_2}(\G_t^\top\G_t).
\]
\noi Hence,
\[
0\preceq S_t^L
\preceq
\|\G_t\|_{\rm sp}^{2}I
\preceq
G^2I .
\]
\noi The EMA recursion $\mathbf{L}_t=(1-\beta_t)\mathbf{L}_{t-1}+\beta_t S_t^L$ is a convex combination. Since $\mathbf{L}_0=\zero$, we have:
\[
0 \preceq \mathbf{L}_t\preceq G^2I .
\]
Thus,
\[
\|\mathbf{L}_t\|_{\rm sp}\le G^2 .
\]
\noi The same argument gives
\[
0\preceq \mathbf{R}_t\preceq G^2I .
\]

\paragraph{Part (d).} Adding the regularization term, $\widetilde{\mathbf{L}}_t=\mathbf{L}_t+\varepsilon I$, and $\widetilde{\mathbf{R}}_t=\mathbf{R}_t+\varepsilon I$. Therefore,
\beq
\varepsilon I  ~~\preceq~~ \widetilde{\mathbf{L}}_t,\widetilde{\mathbf{R}}_t
~~\preceq~~ (\varepsilon+G^2)I. \nn
\eeq

\paragraph{Part (e).} For any matrix $\M_t$, $\langle \M_t,\mathcal Q_t(\M_t)\rangle=\|
\widetilde{\mathbf{L}}_t^{-p/2}\M_t
\widetilde{\mathbf{R}}_t^{-q/2}
\|_{\fro}^{2}$. Because $\widetilde L_t\preceq\varpi I$, and $\widetilde R_t\preceq\varpi I$, we have
\beq
\lambda_{\min} (\widetilde L_t^{-p/2}) \ge \varpi^{-p/2},~~ \lambda_{\min}
(\widetilde R_t^{-q/2})
\ge
\varpi^{-q/2}.\nn
\eeq
This leads to:
\beq
\langle \M_t,\mathcal Q_t(\M_t)\rangle~~\ge~~
\varpi^{-p-q} \|\M_t\|_{\fro}^{2} ~~=~~ \frac1h
\|\M_t\|_{\fro}^{2},\nn
\eeq
where $h=\varpi^{p+q}=\sqrt{\varpi}$. Similarly,
\beq
\|\mathcal Q_t(\M_t)\|_{\fro}~~\le~~
\|\widetilde L_t^{-p}\|_{\rm sp}
\|\M_t\|_{\fro}
\|\widetilde R_t^{-q}\|_{\rm sp}
~~\le~~ \varepsilon^{-p-q}
\|\M_t\|_{\fro}
~~=~~ \varepsilon^{-1/2}\|\M_t\|_{\fro}.\nn
\eeq

\paragraph{Part (f).} We define $C_a:=G^2\varepsilon^{-3/2}$. From the EMA recursion, $\mathbf{L}_t-\mathbf{L}_{t-1} = 
\beta_t(S_t^L- \mathbf{L}_{t-1})$, and therefore
\beq
\|\mathbf{L}_t-\mathbf{L}_{t-1}\|_{\rm sp} \le 2\beta_tG^2. \nn
\eeq
Likewise,
\beq
\|\mathbf{R}_t-\mathbf{R}_{t-1}\|_{\rm sp}
\le
2\beta_tG^2 .\nn
\eeq
\noi Using the fractional-power perturbation lemma,
\beq
\| \widetilde{\mathbf{L}}_t^{-p} - \widetilde{\mathbf{L}}_{t-1}^{-p}
\|_{\rm sp} \le 2p\beta_tG^2\varepsilon^{-p-1},~~\|
\widetilde{\mathbf{R}}_t^{-q} - \widetilde{\mathbf{R}}_{t-1}^{-q}
\|_{\rm sp} \le 2q\beta_tG^2\varepsilon^{-q-1}.\nn
\eeq
Therefore, we have:
\beq
\|\mathcal Q_t(\M_t)-\mathcal Q_{t-1}(\M_t)\|_{\fro}
~~\le~~ 2p\beta_tG^2
\varepsilon^{-p-q-1}\|\M_t\|_{\fro} + 2q\beta_tG^2
\varepsilon^{-p-q-1}\|\M_t\|_{\fro}
~~=~~ G^2\varepsilon^{-3/2} \beta_t\|\M_t\|_{\fro}.\nn
\eeq

\end{proof}

\subsection{Proof of Lemma \ref{lemma:one:step:descent}}
\label{app:lemma:one:step:descent}

\begin{proof}

Recall that
\[
\kappa:=1+\widehat G^2,\qquad
\varpi:=\varepsilon+G^2,\qquad
h:=\varpi^{p+q},\qquad
C_a:=G^2\varepsilon^{-3/2}.
\]

By Lemma~\ref{lemma:preliminary}(f),
\[
\|\mathcal Q_t(\M_t)-\mathcal Q_{t-1}(\M_t)\|_{\fro}
\le
C_a\beta_t\|\M_t\|_{\fro}.
\]

Moreover, by Assumption~\ref{ass:sigma}, Jensen's inequality gives
\[
\|\nabla f(\X_t)\|_{\fro}
=
\|\mathbb E_t[\G_t]\|_{\fro}
\le
\mathbb E_t\|\G_t\|_{\fro}
\le
\hat G .
\]

Define the descent component
\[
\Delta_t
:=
-\theta
\left\langle
\nabla f(\X_t),
\mathcal Q_t(\M_t)
\right\rangle .
\]

Since \(\mathcal Q_{t-1}\) is
\(\mathcal F_{t-1}\)-measurable, it can be treated as a fixed
linear operator when taking conditional expectation with respect to
\(\mathcal F_{t-1}\).

We first separate the variation of the time-varying preconditioner:
\[
\begin{aligned}
\Delta_t
&=
-\theta
\left\langle
\nabla f(\X_t),
(\mathcal Q_t-\mathcal Q_{t-1})(\M_t)
\right\rangle
\\
&\quad
-\theta
\left\langle
\nabla f(\X_t),
\mathcal Q_{t-1}(\M_t)
\right\rangle
\\
&=:A_t+B_t .
\end{aligned}
\]

For the first term, using Lemma~\ref{lemma:preliminary}(f) and
Young's inequality, we obtain
\[
\begin{aligned}
A_t
&\le
\theta C_a\beta_t
\|\nabla f(\X_t)\|_{\fro}
\|\M_t\|_{\fro}
\\
&\le
\frac{\theta\alpha}{4h}
\|\nabla f(\X_t)\|_{\fro}^{2}
+
\frac{\theta h C_a^2}{\alpha}
\beta_t^2
\|\M_t\|_{\fro}^{2}.
\end{aligned}
\]

The second term above is exactly the second component of \(E_t\).

Next, using
\[
\M_t=(1-\alpha)\M_{t-1}+\alpha\G_t ,
\]
we have
\[
\begin{aligned}
B_t
&=
-\theta(1-\alpha)
\left\langle
\nabla f(\X_t),
\mathcal Q_{t-1}(\M_{t-1})
\right\rangle
\\
&\quad
-\theta\alpha
\left\langle
\nabla f(\X_t),
\mathcal Q_{t-1}(\G_t)
\right\rangle .
\end{aligned}
\]

For the stochastic-gradient term, since
\(\mathcal Q_{t-1}\) is
\(\mathcal F_{t-1}\)-measurable,
\[
\begin{aligned}
&
\mathbb E_t
\left\langle
\nabla f(\X_t),
\mathcal Q_{t-1}(\G_t)
\right\rangle
\\
&=
\left\langle
\nabla f(\X_t),
\mathcal Q_{t-1}
(\mathbb E_t[\G_t])
\right\rangle
\\
&=
\left\langle
\nabla f(\X_t),
\mathcal Q_{t-1}
(\nabla f(\X_t))
\right\rangle
\\
&\ge
\frac1h
\|\nabla f(\X_t)\|_{\fro}^{2},
\end{aligned}
\]
where Lemma~\ref{lemma:preliminary}(e) is applied.

It remains to handle the momentum propagation term. Adding and
subtracting \(\nabla f(\X_{t-1})\), we obtain
\[
\begin{aligned}
&
-\theta(1-\alpha)
\left\langle
\nabla f(\X_t),
\mathcal Q_{t-1}(\M_{t-1})
\right\rangle
\\
&=
(1-\alpha)\Delta_{t-1}
\\
&\quad+
\theta(1-\alpha)
\left\langle
\nabla f(\X_{t-1})
-\nabla f(\X_t),
\mathcal Q_{t-1}(\M_{t-1})
\right\rangle .
\end{aligned}
\]

By the update rule,
\[
\X_t-\X_{t-1}
=
-\theta\gamma_{t-1}
\mathcal Q_{t-1}(\M_{t-1}),
\]
and since \(\gamma_t\le\gamma_{t-1}\), smoothness gives
\[
\begin{aligned}
&
\theta(1-\alpha)\gamma_t
\left\langle
\nabla f(\X_{t-1})
-\nabla f(\X_t),
\mathcal Q_{t-1}(\M_{t-1})
\right\rangle
\\
&\le
\theta\gamma_{t-1}
\|\nabla f(\X_{t-1})
-\nabla f(\X_t)\|_{\fro}
\|\mathcal Q_{t-1}(\M_{t-1})\|_{\fro}
\\
&\le
L\theta^2\gamma_{t-1}^{2}
\|\mathcal Q_{t-1}(\M_{t-1})\|_{\fro}^{2}
\\
&\le
\frac{L\theta^2}{\varepsilon}
\gamma_{t-1}^{2}
\|\M_{t-1}\|_{\fro}^{2}.
\end{aligned}
\]

Therefore,
\[
\mathbb E[\gamma_t\Delta_t]
\le
(1-\alpha)
\mathbb E[\gamma_t\Delta_{t-1}]
+
\mathbb E[E_t-D_t].
\]

Furthermore, using
\[
\|\nabla f(\X_{t-1})\|_{\fro}\le \hat G
\]
and Lemma~\ref{lemma:preliminary}(e),
\[
\begin{aligned}
|\Delta_{t-1}|
&\le
\theta
\|\nabla f(\X_{t-1})\|_{\fro}
\|\mathcal Q_{t-1}(\M_{t-1})\|_{\fro}
\\
&\le
\theta\hat G
\varepsilon^{-1/2}
\|\M_{t-1}\|_{\fro}
\\
&\le
\frac{\theta\hat G^2}{\sqrt{\varepsilon}}
=:c_0 .
\end{aligned}
\]

Since \(\gamma_t\le\gamma_{t-1}\),
\[
\gamma_t\Delta_{t-1}
\le
\gamma_{t-1}\Delta_{t-1}
+
c_0(\gamma_{t-1}-\gamma_t).
\]

Applying the above inequality recursively yields
\[
\begin{aligned}
\mathbb E[\gamma_t\Delta_t]
\le&
\mathbb E
\sum_{j=1}^{t}
(1-\alpha)^{t-j}
(E_j-D_j)
\\
&+
(1-\alpha)c_0
\mathbb E
\sum_{j=2}^{t}
(1-\alpha)^{t-j}
(\gamma_{j-1}-\gamma_j).
\end{aligned}
\]

Finally, by \(L\)-smoothness,
\[
\begin{aligned}
f(\X_{t+1})-f(\X_t)
&\le
\gamma_t\Delta_t
+
\frac L2
\theta^2\gamma_t^2
\|\mathcal Q_t(\M_t)\|_{\fro}^{2}
\\
&\le
\gamma_t\Delta_t
+
\frac{L\theta^2}{2\varepsilon}
\gamma_t^2
\|\M_t\|_{\fro}^{2}.
\end{aligned}
\]

Taking expectations and using
\[
c_3=\frac{L\theta^2}{2\varepsilon}
\]
completes the proof.

\end{proof}

\subsection{Proof of Lemma \ref{lemma:log:bound}}
\label{app:lemma:log:bound}

\begin{proof}

We first define
\[
\kappa:=1+\widehat G^2 .
\]

For convenience, we recall that
\[
\gamma_t=(1+\MMM_{t-1})^{-1/2},
\qquad
\MMM_t=\sum_{i=1}^{t}\|\M_i\|_{\fro}^{2},
\qquad
\GGG_t=\sum_{i=1}^{t}\|\G_i\|_{\fro}^{2}.
\]

Since $\MMM_T\le \GGG_T$ by Lemma~\ref{lemma:preliminary}(a), and
$\|\G_t\|_{\fro}\le \widehat G$, we have
\[
\MMM_T\le T\widehat G^2 .
\]

\paragraph{Part (a).}

By the definition of $\gamma_t$,
\[
\sum_{t=1}^{T}\gamma_t^2\|\M_t\|_{\fro}^{2}
=
\sum_{t=1}^{T}
\frac{\|\M_t\|_{\fro}^{2}}
{1+\MMM_{t-1}} .
\]

Applying Lemma~\ref{lemma:ema-denominator}, we obtain
\[
\sum_{t=1}^{T}\gamma_t^2\|\M_t\|_{\fro}^{2}
\le
\kappa\ln(1+\MMM_T).
\]

Using $\MMM_T\le T\widehat G^2$, we have
\[
\sum_{t=1}^{T}\gamma_t^2\|\M_t\|_{\fro}^{2}
\le
\kappa\ln(1+T\widehat G^2).
\]

Since $\widehat G$ is independent of $T$,
\[
\ln(1+T\widehat G^2)
\le
\ln(e+\widehat G^2)+\ln(e+T)
\le
\ln(e+\widehat G^2)\ln(e+T),
\]
where the last inequality follows from $\ln(e+T)\ge1$.
Therefore,
\[
\sum_{t=1}^{T}\gamma_t^2\|\M_t\|_{\fro}^{2}
\le
(1+\widehat G^2)\ln(e+\widehat G^2)\ln(e+T),
\]
which proves
\[
\sum_{t=1}^{T}\gamma_t^2\|\M_t\|_{\fro}^{2}
\le C_b\ln(e+T).
\]

\paragraph{Part (b).}

Recall that
\[
E_t
=
c_1\gamma_{t-1}^{2}\|\M_{t-1}\|_{\fro}^{2}
+
c_2\beta_t^2\gamma_t\|\M_t\|_{\fro}^{2}.
\]

For any nonnegative sequence $\{a_t\}$,
\[
\sum_{t=1}^{T}\sum_{j=1}^{t}
(1-\alpha)^{t-j}a_j
\le
\frac1\alpha\sum_{j=1}^{T}a_j .
\]

Hence,
\[
\sum_{t=1}^{T}\sum_{j=1}^{t}
(1-\alpha)^{t-j}E_j
\le
\frac1\alpha
\sum_{t=1}^{T}E_t .
\]

We first consider the first component of $E_t$.
By Part (a),
\[
\sum_{t=1}^{T}
\gamma_{t-1}^{2}\|\M_{t-1}\|_{\fro}^{2}
\le
\kappa\ln(e+\GGG_T).
\]

Therefore,
\[
c_1
\sum_{t=1}^{T}
\gamma_{t-1}^{2}\|\M_{t-1}\|_{\fro}^{2}
\le
c_1\kappa\ln(e+\GGG_T).
\]

For the second component, Lemma~\ref{lemma:closed-loop-coefficient}
gives
\[
\beta_t^2
\le
\frac{\kappa}{1+\GGG_{t-1}} .
\]

Since $\gamma_t\le1$, we have
\[
\beta_t^2\gamma_t\|\M_t\|_{\fro}^{2}
\le
\kappa
\frac{\|\M_t\|_{\fro}^{2}}
{1+\GGG_{t-1}} .
\]

It remains to bound the last summation.
Using the EMA recursion,
\[
\M_t
=
\alpha\sum_{i=1}^{t}
(1-\alpha)^{t-i}\G_i ,
\]
and the convexity of the squared Frobenius norm,
\[
\|\M_t\|_{\fro}^{2}
\le
\alpha
\sum_{i=1}^{t}
(1-\alpha)^{t-i}
\|\G_i\|_{\fro}^{2}.
\]

Consequently,
\[
\begin{aligned}
\sum_{t=1}^{T}
\frac{\|\M_t\|_{\fro}^{2}}
{1+\GGG_{t-1}}
&\le
\alpha
\sum_{i=1}^{T}
\|\G_i\|_{\fro}^{2}
\sum_{t=i}^{T}
\frac{(1-\alpha)^{t-i}}
{1+\GGG_{t-1}}
\\
&\le
\alpha
\sum_{i=1}^{T}
\|\G_i\|_{\fro}^{2}
\sum_{t=i}^{T}
(1-\alpha)^{t-i}
\\
&\le
\sum_{i=1}^{T}
\|\G_i\|_{\fro}^{2}
\\
&=
\GGG_T .
\end{aligned}
\]

The above bound is not sufficient for logarithmic control.
Instead, using the monotonicity of $\GGG_t$ and
$\|\G_t\|_{\fro}\le\widehat G$, we have

\[
1+\GGG_{t-1}
\ge
\frac{1+\GGG_t}{1+\widehat G^2}.
\]

Therefore,

\[
\frac1{1+\GGG_{t-1}}
\le
\kappa
\frac1{1+\GGG_t},
\]

where $\kappa=1+\widehat G^2$.
Hence,

\[
\begin{aligned}
\sum_{t=1}^{T}
\frac{\|\M_t\|_{\fro}^{2}}
{1+\GGG_{t-1}}
&\le
\kappa
\sum_{t=1}^{T}
\frac{\|\M_t\|_{\fro}^{2}}
{1+\GGG_t}.
\end{aligned}
\]

Applying Lemma~\ref{lemma:ema-denominator}
with the denominator $1+\GGG_t$ yields

\[
\sum_{t=1}^{T}
\frac{\|\M_t\|_{\fro}^{2}}
{1+\GGG_{t-1}}
\le
\kappa^2
\ln(1+\GGG_T).
\]

Consequently,

\[
\sum_{t=1}^{T}
\beta_t^2\gamma_t\|\M_t\|_{\fro}^{2}
\le
\kappa^3
\ln(1+\GGG_T).
\]

Combining the two components gives

\[
\sum_{t=1}^{T}
\sum_{j=1}^{t}
(1-\alpha)^{t-j}E_j
\le
\left(
\frac{c_1\kappa}{\alpha}
+
\frac{c_2\kappa^3}{\alpha}
\right)
\ln(1+\GGG_T).
\]

Including the smoothness contribution,
\[
c_3
\sum_{t=1}^{T}
\gamma_t^2\|\M_t\|_{\fro}^{2}
\le
c_3C_b\ln(e+T),
\]

we finally obtain

\[
\sum_{t=1}^{T}
\sum_{j=1}^{t}
(1-\alpha)^{t-j}E_j
\le
C_e\ln(1+T\widehat G^2),
\]

where

\[
C_e
=
\frac{c_1\kappa}{\alpha}
+
\frac{c_2\kappa^3}{\alpha}
+
c_3\kappa .
\]

\paragraph{Part (c).}

Summing the one-step descent inequality in
Lemma~\ref{lemma:one:step:descent} over $t=1,\ldots,T$ gives
\[
\begin{aligned}
c_4\alpha
\mathbb E[H_T]
\le&
f(\X_1)-f_*
\\
&+
\mathbb E
\sum_{t=1}^{T}
c_3\gamma_t^2\|\M_t\|_{\fro}^{2}
\\
&+
\mathbb E
\sum_{t=1}^{T}
\sum_{j=1}^{t}
(1-\alpha)^{t-j}E_j
\\
&+
(1-\alpha)c_0
\sum_{t=1}^{T}
\sum_{j=2}^{t}
(1-\alpha)^{t-j}
(\gamma_{j-1}-\gamma_j).
\end{aligned}
\]

The last term is bounded by
\[
\sum_{t=1}^{T}
\sum_{j=2}^{t}
(1-\alpha)^{t-j}
(\gamma_{j-1}-\gamma_j)
\le
\frac1\alpha .
\]

Using Parts (a) and (b), we obtain
\[
\begin{aligned}
c_4\alpha
\mathbb E[H_T]
\le&
f(\X_1)-f_*
+
C_e\ln(1+T\widehat G^2)
\\
&+
c_3C_b\ln(e+T)
+
\frac{(1-\alpha)c_0}{\alpha}.
\end{aligned}
\]

Since $\widehat G$ is independent of $T$, the logarithmic terms can be
absorbed into $\ln(e+T)$ by enlarging the constant. Hence,
\[
\mathbb E[H_T]
\le
C_h\ln(e+T),
\]
where
\[
C_h
=
\frac{
f(\X_1)-f_*
+C_e
+c_3C_b
+\frac{(1-\alpha)c_0}{\alpha}
}
{c_4\alpha}.
\]

This completes the proof.
\end{proof}

\subsection{Proof of Lemma \ref{lemma:Weighted:energy:bridge}}
\label{app:lemma:Weighted:energy:bridge}

\begin{proof}

Recall that
\[
H_T
=
\sum_{t=1}^{T}
\gamma_t\|\nabla f(\X_t)\|_{\fro}^{2},
\]
and
\[
\gamma_t=(1+\MMM_{t-1})^{-1/2}.
\]

Since $\MMM_t$ is nondecreasing, for every $t\le T$,
\[
\MMM_{t-1}\le \MMM_T .
\]
Therefore,
\[
\gamma_t
=
\frac{1}{\sqrt{1+\MMM_{t-1}}}
\ge
\frac{1}{\sqrt{1+\MMM_T}} .
\]

Consequently,
\[
\begin{aligned}
H_T
&=
\sum_{t=1}^{T}
\gamma_t
\|\nabla f(\X_t)\|_{\fro}^{2}
\\
&\ge
\frac{1}{\sqrt{1+\MMM_T}}
\sum_{t=1}^{T}
\|\nabla f(\X_t)\|_{\fro}^{2}
\\
&=
\frac{A_T}{\sqrt{1+\MMM_T}} .
\end{aligned}
\]

Hence,
\[
A_T
\le
H_T\sqrt{1+\MMM_T},
\]
and taking square roots gives
\[
\sqrt{A_T}
\le
\sqrt{H_T}
(1+\MMM_T)^{1/4}.
\]

Taking expectation on both sides and applying
Cauchy--Schwarz inequality yields
\[
\begin{aligned}
\mathbb E\sqrt{A_T}
&\le
\mathbb E
\left[
\sqrt{H_T}
(1+\MMM_T)^{1/4}
\right]
\\
&\le
(\mathbb E H_T)^{1/2}
\left(
\mathbb E\sqrt{1+\MMM_T}
\right)^{1/2}.
\end{aligned}
\]

By Lemma~\ref{lemma:log:bound}(c),
\[
\mathbb E H_T
\le
C_h\ln(e+T).
\]

Moreover, Lemma~\ref{lemma:preliminary}(a) gives
\[
\MMM_T\le \GGG_T .
\]
Therefore,
\[
\mathbb E\sqrt{1+\MMM_T}
\le
\mathbb E\sqrt{1+\GGG_T}
=
Z_T .
\]

Combining the above inequalities, we obtain
\[
\begin{aligned}
\mathbb E\sqrt{A_T}
&\le
(C_h\ln(e+T))^{1/2}
Z_T^{1/2}
\\
&=
C_h^{1/2}
Z_T^{1/2}
\ln^{1/2}(e+T).
\end{aligned}
\]

Letting
\[
C_Z:=C_h^{1/2},
\]
which is independent of $T$ and $\sigma$, we conclude that
\[
\mathbb E\sqrt{A_T}
\le
C_Z
Z_T^{1/2}
\ln^{1/2}(e+T).
\]

This completes the proof.
\end{proof}

\subsection{Proof of Theorem \ref{the:1}}
\label{app:the:1}

\begin{proof}

By Lemma~\ref{lemma:Weighted:energy:bridge},
we have
\[
\mathbb E\sqrt{A_T}
\le
C_Z Z_T^{1/2}\ln^{1/2}(e+T),
\]
where
\[
C_Z=\sqrt{C_h},
\]
and
\[
C_h=
\frac{
f(\X_1)-f_*
+C_e
+c_3C_b
+\frac{(1-\alpha)c_0}{\alpha}
}
{c_4\alpha}.
\]

Recall that
\[
Z_T=\mathbb E\sqrt{1+\GGG_T},
\]
where
\[
\GGG_T=\sum_{t=1}^{T}\|\G_t\|_{\fro}^{2}.
\]

Let
\[
\varepsilons_t
=
\G_t-\nabla f(\X_t).
\]
Then
\[
\G_t=\nabla f(\X_t)+\varepsilons_t.
\]

Using
\[
\|a+b\|_{\fro}^{2}
\le
2\|a\|_{\fro}^{2}
+
2\|b\|_{\fro}^{2},
\]
we obtain
\[
\begin{aligned}
1+\GGG_T
&=
1+
\sum_{t=1}^{T}
\|\nabla f(\X_t)+\varepsilons_t\|_{\fro}^{2}
\\
&\le
1+
2\sum_{t=1}^{T}
\|\nabla f(\X_t)\|_{\fro}^{2}
+
2\sum_{t=1}^{T}
\|\varepsilons_t\|_{\fro}^{2}
\\
&=
1+2A_T
+
2\sum_{t=1}^{T}
\|\varepsilons_t\|_{\fro}^{2}.
\end{aligned}
\]

Therefore,
\[
\sqrt{1+\GGG_T}
\le
1+\sqrt2\sqrt{A_T}
+
\sqrt2
\sqrt{
\sum_{t=1}^{T}
\|\varepsilons_t\|_{\fro}^{2}
}.
\]

Taking expectations gives
\[
Z_T
\le
1
+
\sqrt2\mathbb E\sqrt{A_T}
+
\sqrt2
\mathbb E
\sqrt{
\sum_{t=1}^{T}
\|\varepsilons_t\|_{\fro}^{2}
}.
\]

For the stochastic noise term, Jensen's inequality and
Assumption~\ref{ass:sigma} yield
\[
\begin{aligned}
\mathbb E
\sqrt{
\sum_{t=1}^{T}
\|\varepsilons_t\|_{\fro}^{2}
}
&\le
\sqrt{
\mathbb E
\sum_{t=1}^{T}
\|\varepsilons_t\|_{\fro}^{2}
}
\\
&\le
\sigma\sqrt T .
\end{aligned}
\]

Hence,
\[
Z_T
\le
1+
\sqrt2C_Z
Z_T^{1/2}\ln^{1/2}(e+T)
+
\sqrt2\sigma\sqrt T .
\]

Define
\[
Y_T:=Z_T^{1/2}.
\]
Then
\[
Y_T^2
-
\sqrt2C_Z\ln^{1/2}(e+T)Y_T
-
(1+\sqrt2\sigma\sqrt T)
\le0 .
\]

Solving the quadratic inequality gives
\[
Y_T
\le
\frac{
\sqrt2C_Z\ln^{1/2}(e+T)
+
\sqrt{
2C_Z^2\ln(e+T)
+
4(1+\sqrt2\sigma\sqrt T)
}
}
{2}.
\]

Therefore,
\[
Z_T^{1/2}
\le
\frac{
\sqrt2C_Z\ln^{1/2}(e+T)
+
\sqrt{
2C_Z^2\ln(e+T)
+
4(1+\sqrt2\sigma\sqrt T)
}
}
{2}.
\]

Substituting this bound into
Lemma~\ref{lemma:Weighted:energy:bridge}, we obtain

\[
\begin{aligned}
\mathbb E\sqrt{A_T}
\le&
\frac{C_Z\ln^{1/2}(e+T)}{2}
\\
&\times
\Bigg(
\sqrt2C_Z\ln^{1/2}(e+T)
\\
&\qquad+
\sqrt{
2C_Z^2\ln(e+T)
+
4(1+\sqrt2\sigma\sqrt T)
}
\Bigg).
\end{aligned}
\]

Now, by Cauchy--Schwarz inequality,

\[
\sum_{t=1}^{T}
\|\nabla f(\X_t)\|_{\fro}
\le
\sqrt{
T
\sum_{t=1}^{T}
\|\nabla f(\X_t)\|_{\fro}^{2}
}
=
\sqrt{TA_T}.
\]

Hence,

\[
\frac1T
\sum_{t=1}^{T}
\|\nabla f(\X_t)\|_{\fro}
\le
\frac{\sqrt{A_T}}{\sqrt T}.
\]

Taking expectations gives

\[
\begin{aligned}
&
\mathbb E
\left[
\frac1T
\sum_{t=1}^{T}
\|\nabla f(\X_t)\|_{\fro}
\right]
\\
&\le
\frac{C_Z\ln^{1/2}(e+T)}
{2\sqrt T}
\\
&\quad\times
\Bigg(
\sqrt2C_Z\ln^{1/2}(e+T)
\\
&\qquad+
\sqrt{
2C_Z^2\ln(e+T)
+
4(1+\sqrt2\sigma\sqrt T)
}
\Bigg).
\end{aligned}
\]

Finally, let
\[
R\sim{\rm Unif}\{1,\ldots,T\}
\]
be independent of the algorithmic randomness. Then

\[
\mathbb E\|\nabla f(\X_R)\|_{\fro}
=
\mathbb E
\left[
\frac1T
\sum_{t=1}^{T}
\|\nabla f(\X_t)\|_{\fro}
\right],
\]

which proves the random-output statement.

When $\sigma=0$, the above inequality reduces to

\[
\begin{aligned}
\mathbb E\|\nabla f(\X_R)\|_{\fro}
\le&
\frac{C_Z\ln^{1/2}(e+T)}
{2\sqrt T}
\\
&\times
\Bigg(
\sqrt2C_Z\ln^{1/2}(e+T)
+
\sqrt{
2C_Z^2\ln(e+T)+4
}
\Bigg),
\end{aligned}
\]

which gives the deterministic-order rate
\[
\widetilde{\mathcal O}(T^{-1/2}).
\]

This completes the proof.
\end{proof}

\section{Implementation Configuration}\label{supp:sec:implementation}
For reproducibility, this section consolidates the implementation choices stated in Section~\ref{sect:experiments}. Matrix parameters use NABLA, while one-dimensional parameters and embeddings use an AdamW fallback. CIFAR-100 uses 20 Newton--Schulz iterations with inverse factors refreshed at every update; the compared realization is NABLA$(64,128)$ against NABLA$(m,n)$. LLaMA-60M compares NABLA$(1,1)$ against NABLA$(m,n)$ under shared momentum construction, inverse-factor postprocessing, update stabilization, and scalar hyperparameters. A release version should additionally report the exact parameter-group rules, Newton--Schulz normalization, damping implementation, random seeds, and software commit.


\begin{figure*}[t!]
\centering
\includegraphics[width=0.96\textwidth]
{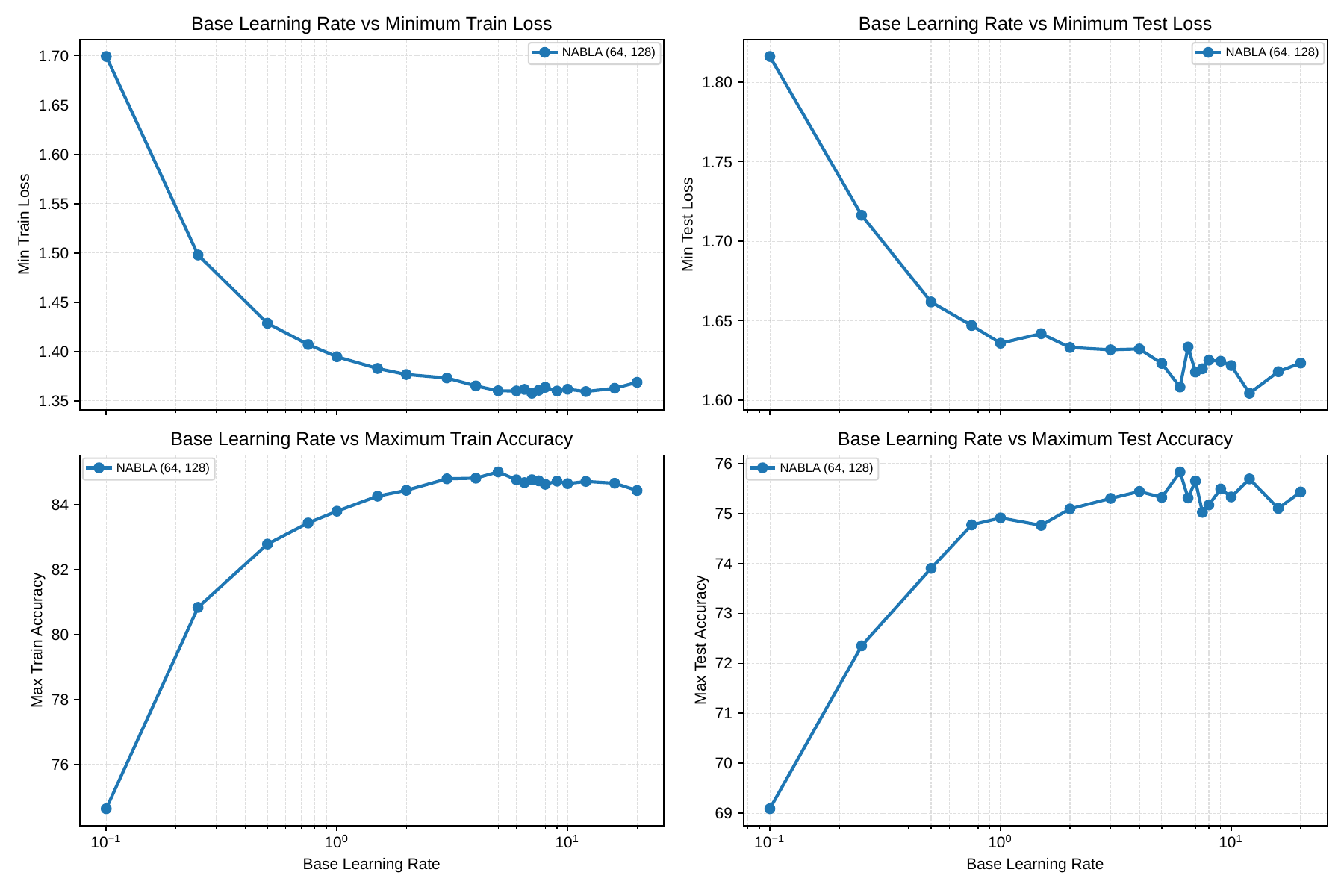}
\caption{Base-learning-rate sensitivity of NABLA$(64,128)$ on CIFAR-100.
The top row reports the minimum training and test losses attained within each
fixed-horizon sweep run, and the bottom row reports the corresponding maximum
training and test accuracies. Lower is better for loss and higher is better
for accuracy. This fixed-seed sweep is used as a sensitivity diagnostic and
is separate from the paired-seed optimizer comparison in
Figure~\ref{fig:training-dynamics}.}
\label{fig:cifar-lr-sensitivity}
\end{figure*}

\end{document}